\documentclass[12pt,oneside]{article}
\usepackage{amsmath,amssymb,amsfonts,amsthm}
\def\goth{\mathfrak}

\newcommand{\T}{{\cal T}}

\newcommand{\Real}{\mathbb R}

\newcommand{\To}{\longrightarrow}

\newcommand{\prof}{\noindent \textit{\textbf{Proof.\:\:}}}

\newcommand{\p}{\pi^{-1}(TM)}
\newcommand {\cp}{\mathfrak{X}(\pi (M))}
\newcommand {\cpp}{\mathfrak{X}(\T M)}

\def\x{{\goth X}(\T M)}

\def\Section#1{\vspace{30truept}\addtocounter{section}{1}\setcounter{thm}{0}
\setcounter{equation}{0}{\noindent\Large\bf
    \arabic{section}.~~#1}\par \vspace{12pt}}

\newtheorem{thm}{Theorem}[section]
\newtheorem{cor}[thm]{Corollary}
\newtheorem{lem}[thm]{Lemma}
\newtheorem{prop}[thm]{Proposition}
\newtheorem{defn}[thm]{Definition}

\newtheorem{rem}[thm]{Remark}

\numberwithin{equation}{section}
\begin{document}
\title{\bf{ CARTAN AND BERWALD CONNECTIONS IN THE PULLBACK FORMALISM} }
\author{{\bf Nabil L. Youssef$^{\dag}$, S. H. Abed$^{\dag}$ and A. Soleiman$^{\ddag}$}}
\date{}
%\thanks{\it Department of Mathematics, etc}
%\pagestyle{fancy}

             % End of preamble and beginning of text.
\maketitle                     % Produces the title.
\vspace{-1.15cm}
\begin{center}
{$^{\dag}$Department of Mathematics, Faculty of Science,\\ Cairo
University, Giza, Egypt}
\end{center}
\vspace{-0.8cm}
\begin{center}
nlyoussef2003@yahoo.fr,\ sabed52@yahoo.fr
\end{center}
\vspace{-0.7cm}
\begin{center}
and
\end{center}
\vspace{-0.7cm}
\begin{center}
{$^{\ddag}$Department of Mathematics, Faculty of Science,\\ Benha
University, Benha,
 Egypt}
\end{center}
\vspace{-0.8cm}
\begin{center}
soleiman@mailer.eun.eg
\end{center}
\smallskip

\vspace{1cm} \maketitle
\smallskip

\noindent{\bf Abstract.} Adopting the pullback approach to global
Finsler geometry, the aim of the present paper is to provide new
intrinsic (coordinate-free) proofs of intrinsic versions of the
existence and uniqueness theorems for the Cartan and Berwald
connections on a Finsler manifold. To accomplish this, the notions
of semispray and nonlinear connection associated with a given
regular connection, in the pullback bundle, is introduced and
investigated. Moreover, it is shown that for the Cartan and Berwald
connections, the associated semispray coincides with the canonical
spray and the associated nonlinear connection coincides with the
Barthel connection. An explicit intrinsic expression relating both
connections is deduced.
\par
Although our treatment is entirely global, the local expressions of
the obtained results, when calculated, coincide with the existing
classical local results.



%%%%%%%%%%%%%%%%%%%%%%%%%%%%%%%%%%%%%%%%%%%%%%%%%%%%%%%% INTRODUCTION %%%%%%%%%%%%%%%%%%%%%%%%%%%%%%%%%%%%%%%%%%%%%%%%%%%%%%%%%%%%%%%%%%%%%%%%%%%

%\vspace{1cm}
%Introduction
\vspace{30truept}\centerline{\Large\bf{Introduction}}\vspace{12pt}
\par
The most well-known and widely used approaches to GLOBAL Finsler
geometry are the Klein-Grifone (KG-) approach (cf. \cite{r21},
\cite{r22}, \cite{r27}) and the pullback (PB-) approach (cf.
\cite{r58}, \cite{r61}, \cite{r74}, \cite{r44}). The universe of the
first approach is the tangent bundle of $\,\T M$ (i.e, $\pi_{\T
M}:T\T M\longrightarrow \T M$), whereas the universe of the second
is the pullback of the tangent bundle $TM$ by $\pi: \T
M\longrightarrow M$  (i.e., $P:\pi^{-1}(TM)\longrightarrow \T M$).
Each of the two approaches has its own geometry which differs
significantly from the geometry of the other (in spite of the
existence of some links between them).
\par
 In Riemannian geometry, there is a canonical linear connection on the
manifold $M$, whereas in Finsler geometry there is a corresponding
canonical linear connection due to E. Cartan. However, this is not a
connection on $M$ but is a connection on $T(\T M)$  (\emph{in the
KG- approach}) or on $\,\pi^{-1}(TM) $ (\emph{in the PB-approach}).
\par
The most important linear connections in Finsler geometry are the
Cartan connection and the Berwald connection. On the other hand,
local Finsler geometry, which is very widespread, is the local
version of the PB-approach. These are among the reasons that
motivated this work. Moreover, to the best of our knowledge there is
no proof, in the PB-approach, of the existence and uniqueness
theorems for the Cartan and Berwald connections from a purely global
perspective.
\par
The  main purpose of the present paper is to provide \textbf{new
intrinsic} (coordinate-free) proofs of intrinsic versions of the
existence and uniqueness theorems for the Cartan and Berwald
connections within the pullback formalism, making simultaneous use
of some concepts and results from the KG-approach. These proofs have
the advantages of being simple, systematic and parallel to and
guided by the Reimannian case. It is worth mentioning here that our
proofs are fundamentally different from that given by P. Dazord
\cite{r61}, which is not purely intrinsic.
\par
The  paper consists of three parts preceded by an introductory
section $(\S 1)$, which provides a brief account of the basic
definitions and concepts necessary for this work. For more details,
we refer to \cite{r44}, \cite{r61}, \cite{r21} and \cite{r22}
\par
In the first part $(\S 2)$, the notions of semispray and nonlinear
connection associated with a given regular connection, in the
pullback bundle, are introduced and investigated.
\par
The second part $(\S 3)$ is devoted to an intrinsic proof of the
existence and uniqueness theorem of the Cartan connection on a
Finsler manifold $(M,L)$ (Theorem \ref{th.1}). For the Cartan
connection, it is shown that the associated semispray coincides with
the canonical spray (Corollary \ref{Cth.1}) and the associated
nonlinear connection coincides with the Barthel connection (Theorem
\ref{c.ba.}). This establishes an important link between the
PB-approach and the KG-approach.
\par
The third and last part $(\S 4)$ provides an intrinsic proof of the
existence and uniqueness theorem of the Berwald connection on
$(M,L)$ (Theorem \ref{bth2.h2}). Moreover, an elegant formula
relating this connection and the Cartan connection is obtained
(Theorem \ref{th.5}). A by-product of the above results is a
characterization of Riemannian and Landsbergian manifolds.
\par
We have to emphasize that without the insertion of the KG-approach
we would have been unable to achieve these results. It should also
be pointed out that the present work is formulated in a prospective
modern coordinate-free form; the local expressions of the obtained
results, when calculated, coincide with the existing classical local
results.
\par
Finally, it is worth noting that there are other connections of
particular importance in Finsler geometry, such as Chern (Rund) and
Hashiguchi connections, which are not treated in the present work.
They merit a separate study that we are currently in the process of
preparing and will be the object of a forthcoming paper.

%%%%%%%%%%%%%%%%%%%%%%%%%%%%%$$$$$$$$$$$$$$$$$$ SECTION 1. Notation and Preliminaries $$$$$$$$$$$$$$$$$$$$$$$$$$$$$$$$$$$$$$$$$$$$$$$$$$$$$$$$$$$$$$$

\Section{Notation and Preliminaries}

In this section, we give a brief account of the basic concepts
 of the pullback formalism necessary for this work. For more
 details, we refer to \cite{r58},\,\cite{r61},\,\cite{r74} and~\,\cite{r44}.
 We make the
assumption that the geometric objects we consider are of class
$C^{\infty}$.\\ The
following notation will be used throughout this paper:\\
 $M$: a real differentiable manifold of finite dimension $n$ and of
class $C^{\infty}$,\\
 $\mathfrak{F}(M)$: the $\Real$-algebra of differentiable functions
on $M$,\\
 $\mathfrak{X}(M)$: the $\mathfrak{F}(M)$-module of vector fields
on $M$,\\
$\pi_{M}:TM\longrightarrow M$: the tangent bundle of $M$,\\
$\pi: \T M\longrightarrow M$: the subbundle of nonzero vectors
tangent to $M$,\\
$V(TM)$: the vertical subbundle of the bundle $TTM$,\\
 $P:\pi^{-1}(TM)\longrightarrow \T M$ : the pullback of the
tangent bundle $TM$ by $\pi$,\\
 $\mathfrak{X}(\pi (M))$: the $\mathfrak{F}(\T M)$-module of
differentiable sections of  $\pi^{-1}(T M)$,\\
$ i_{X}$ : the interior product with respect to  $X
\in\mathfrak{X}(M)$,\\
$df$ : the exterior derivative  of $f$,\\
$ d_{L}:=[i_{L},d]$, $i_{L}$ being the interior derivative with
respect to a vector form $L$.

\par Elements  of  $\mathfrak{X}(\pi (M))$ will be called
$\pi$-vector fields and will be denoted by barred letters
$\overline{X} $. Tensor fields on $\pi^{-1}(TM)$ will be called
$\pi$-tensor fields. The fundamental $\pi$-vector field is the
$\pi$-vector field $\overline{\eta}$ defined by
$\overline{\eta}(u)=(u,u)$ for all $u\in \T M$.

We have the following short exact sequence of vector bundles,
relating the tangent bundle $T(\T M)$ and the pullback bundle
$\pi^{-1}(TM)$:\vspace{-0.1cm}
$$0\longrightarrow
 \pi^{-1}(TM)\stackrel{\gamma}\longrightarrow T(\T M)\stackrel{\rho}\longrightarrow
\pi^{-1}(TM)\longrightarrow 0 ,\vspace{-0.1cm}$$
 where the bundle morphisms $\rho$ and $\gamma$ are defined respectively by
$\rho := (\pi_{\T M},d\pi)$ and $\gamma (u,v):=j_{u}(v)$, where
$j_{u}$  is the natural isomorphism $j_{u}:T_{\pi_{M}(v)}M
\longrightarrow T_{u}(T_{\pi_{M}(v)}M)$. The vector $1$-form $J$ on
$TM$ defined by $J:=\gamma\circ\rho$ is called the natural almost
tangent structure of $T M$. Clearly,
$\textit{Im}\,J=\textit{Ker}\,J=V(TM)$. The vertical vector field
$\mathcal{C}$ on $TM$ defined by
$\mathcal{C}:=\gamma\circ\overline{\eta} $ is called the fundamental
or the canonical (Liouville) vector field.

Let $\nabla$ be  a linear connection (or simply a connection) on the
pullback bundle $\pi^{-1}(TM)$.
 We associate with
$\nabla$ the map \vspace{-0.1cm}
$$K:T \T M\longrightarrow \pi^{-1}(TM):X\longmapsto \nabla_X \overline{\eta}
,\vspace{-0.1cm}$$ called the connection (or the deflection) map of
$\nabla$. A tangent vector $X\in T_u (\T M)$ is said to be
horizontal if $K(X)=0$ . The vector space $H_u (\T M)= \{ X \in T_u
(\T M) : K(X)=0 \}$ of the horizontal vectors
 at $u \in  \T M$ is called the horizontal space to $M$ at $u$  .
   The connection $\nabla$ is said to be regular if
\begin{equation}\label{direct sum}
T_u (\T M)=V_u (\T M)\oplus H_u (\T M) \qquad \forall u\in \T M .
\end{equation}
\par If $M$ is endowed with a regular connection, then the vector bundle
   maps
\begin{eqnarray*}
 \gamma &:& \pi^{-1}(T M)  \To V(\T M), \\
   \rho |_{H(\T M)}&:&H(\T M) \To \pi^{-1}(TM), \\
   K |_{V(\T M)}&:&V(\T M) \To \pi^{-1}(T M)
\end{eqnarray*}
 are vector bundle isomorphisms.
   Let us denote
 $\beta:=(\rho |_{H(\T M)})^{-1}$,
then \vspace{-0.2cm}
   \begin{align}\label{fh1}
    \rho\circ\beta = id_{\pi^{-1} (TM)}, \quad  \quad
       \beta\circ\rho =\left\{
                                \begin{array}{ll}
                                          id_{H(\T M)} & {\,\, on\,\,   H(\T M)} \\
                                         0 & {\,\, on \,\,   V(\T M)}
                                       \end{array}
                                     \right.\vspace{-0.2cm}
\end{align}
The map $\beta$ will be called the horizontal map of the connection
$D$.
\par According to the direct sum decomposition (\ref{direct
sum}), a regular connection $\nabla$ induces a horizontal projector
$h_{\nabla}$ and a vertical projector $v_{\nabla}$, given by
\begin{equation}\label{proj.}
h_{\nabla}=\beta\circ\rho ,  \ \ \ \ \ \ \ \ \ \ \
v_{\nabla}=I-\beta\circ\rho,
\end{equation}
where $I$ is the identity endomorphism on $T(TM)$: $I=id_{T(TM)}$.
\par
 The (classical)  torsion tensor $\textbf{T}$  of the connection
$\nabla$ is defined by
$$\textbf{T}(X,Y)=\nabla_X \rho Y-\nabla_Y\rho X -\rho [X,Y] \quad
\forall\,X,Y\in \mathfrak{X} (\T M).$$ The horizontal ((h)h-) and
mixed ((h)hv-) torsion tensors, denoted by $Q $ and $ T $
respectively, are defined by \vspace{-0.2cm}
$$Q (\overline{X},\overline{Y})=\textbf{T}(\beta \overline{X}\beta \overline{Y}),
\, \,\, T(\overline{X},\overline{Y})=\textbf{T}(\gamma
\overline{X},\beta \overline{Y}) \quad \forall \,
\overline{X},\overline{Y}\in\mathfrak{X} (\pi (M)).\vspace{-0.2cm}$$
\par The (classical) curvature tensor  $\textbf{K}$ of the connection
$\nabla$ is defined by
 $$ \textbf{K}(X,Y)\rho Z=-\nabla_X \nabla_Y \rho Z+\nabla_Y \nabla_X \rho Z+\nabla_{[X,Y]}\rho Z
  \quad \forall\, X,Y, Z \in \mathfrak{X} (\T M).$$
The horizontal (h-), mixed (hv-) and vertical (v-) curvature
tensors, denoted by $R$, $P$ and $S$ respectively, are defined by
$$R(\overline{X},\overline{Y})\overline{Z}=\textbf{K}(\beta
\overline{X}\beta \overline{Y})\overline{Z},\quad
P(\overline{X},\overline{Y})\overline{Z}=\textbf{K}(\beta
\overline{X},\gamma \overline{Y})\overline{Z},\quad
S(\overline{X},\overline{Y})\overline{Z}=\textbf{K}(\gamma
\overline{X},\gamma \overline{Y})\overline{Z}.$$

The contracted curvature tensors, denoted by $\widehat{R}$,
$\widehat{P}$ and $\widehat{S}$ respectively, are also known as the
 (v)h-, (v)hv- and (v)v-torsion tensors and are defined by
$$\widehat{R}(\overline{X},\overline{Y})={R}(\overline{X},\overline{Y})\overline{\eta},\quad
\widehat{P}(\overline{X},\overline{Y})={P}(\overline{X},\overline{Y})\overline{\eta},\quad
\widehat{S}(\overline{X},\overline{Y})={S}(\overline{X},\overline{Y})\overline{\eta}.$$

\par Let $(M,L)$ be a Finsler manifold and $g$ the associated Finsler metric.
 For a regular connection $\nabla$,  we
define
$$R(\overline{X},\overline{Y},\overline{Z}, \overline{W}):
=g(R(\overline{X},\overline{Y})\overline{Z}, \overline{W}),\,
\cdots, \, S(\overline{X},\overline{Y},\overline{Z}, \overline{W}):
=g(S(\overline{X},\overline{Y})\overline{Z}, \overline{W}).$$
\par We terminate this section
by some concepts and results from the Klein-Grifone formalism. For
more details, we refer to \cite{r21}, \cite{r22}, \cite{r27} and 
\cite{n5}.

 A semispray  on $M$ is a vector field $X$ on $TM$,
 $C^{\infty}$ on $\T M$, $C^{1}$ on $TM$, such that
$\rho\circ X = \overline{\eta}$. A semispray $X$ which is
homogeneous of degree $2$ in the directional argument
($[\mathcal{C},X]= X $) is called a spray.

\begin{prop}{\em{\cite{r27}}}\label{spray} Let $(M,L)$ be a Finsler manifold. The vector field $G$
 defined by $i_{G}\Omega =-dE$ is a spray, where
 $E:=\frac{1}{2}L^{2}$ is the energy function and $\Omega:=dd_{J}E$.
 Such a spray is called the canonical spray.
 \end{prop}

A nonlinear connection on $M$ is a vector $1$-form $\Gamma$ on $TM$,
$C^{\infty}$ on $\T M$, $C^{0}$ on $TM$, such that
$$J \Gamma=J, \quad\quad \Gamma J=-J .$$
The horizontal and vertical projectors $h_{\Gamma}$\,  and
$v_{\Gamma}$ associated with $\Gamma$ are defined by
   $h_{\Gamma}:=\frac{1}{2} (I+\Gamma),\, v_{\Gamma}:=\frac{1}{2}
 (I-\Gamma).$
Thus $\Gamma$ gives rise to the direct sum decomposition $T\T M=
H(\T M)\oplus V(\T M)$, where $H(\T M):=Im \, h_{\Gamma} = Ker\,
v_{\Gamma} $, $V(\T M):= Im \, v_{\Gamma}=Ker \, h_{\Gamma}$. We
have $ J\circ h_{\Gamma} =J, \, h_{\Gamma}\circ J=0, \, J\circ
v_{\Gamma}=0, \, v_{\Gamma}\circ J=J.$ A nonlinear connection
$\Gamma$ is homogeneous if $[\mathcal{C},\Gamma]=0$. The torsion $t$
of a nonlinear connection $\Gamma$ is the vector $2$-form  on $TM$
defined by $t:=\frac{1}{2} [J,\Gamma]$. A nonlinear connection
$\Gamma$ is said to be conservative if $d_{h_{\Gamma}}\,E=0$. With
any given nonlinear connection $\Gamma$, one can associate a
semispray $S$ which is horizontal with respect to $\Gamma$, namely,
$S=h_{\Gamma}S'$, where $S'$ is an arbitrary semispray. Moreover, if
$\Gamma$ is homogeneous, then its associated semispray is a spray.

\begin{thm} \label{th.9a} {\em{\cite{r22}}} On a Finsler manifold $(M,L)$, there exists a unique
conservative homogenous nonlinear  connection  with zero torsion. It
is given by\,{\em:} \vspace{-0.3cm} $$\Gamma =
[J,G],\vspace{-0.3cm}$$ where $G$ is the canonical spray.\\
 Such a nonlinear connection is called the canonical connection, or the Barthel connection, associated with $(M,L)$.
\end{thm}
\par
It should be noted that the semispray associated with the Barthel
connection is a spray, which is the canonical spray.

%%%%%%%%%%%%%%%%%%%%%%%%%%%%%%%%%%%%%% SECTION 2. Regular connection on the Pullback Bundle %%%%%%%%%%%%%%%%%%%%%%%%%%%%%%%%%%%%%%%%%%%%%%%%%

\Section{Regular Connections in the Pullback Bundle}

In this section, the semispray and the nonlinear connection
associated with a given regular connection on $\pi^{-1}(TM)$  are
introduced and investigated. \vspace{5pt}
\par
 The following lemma is useful for subsequent use.\vspace{-0.2cm}
\begin{lem}\label{iso}Let $D$ be a regular connection on $\pi^{-1}(TM)$ with horizontal map $\beta$.
Let $S$ be an arbitrary  semispray on $M$. Then, we have
\begin{description}
    \item[(a)] $\rho\, X= \rho\,[J X,S]$, for every $X\in\cpp$,
    \item[(b)] $\overline{X}= \rho\,[\gamma \overline{X},S]$, for
    every $\overline{X}\in\cp$.
\end{description}
\end{lem}

%\begin{proof}
\prof It is known that \cite{r21} any vertical vector field
    $JX$ can be written in the form  $J X= J[J X,S]$,  where
     $S$ is an arbitrary  semispray.\newline
\textbf{(a)} As $J=\gamma\circ\rho\,$  and $\,\gamma: \p
\longrightarrow V(\T M)$ is an isomorphism, then \textbf{(a)}
follows.\newline \textbf{(b)} Follows from \textbf{(a)} by setting
$X=\beta \overline{X}$ and noting that $\rho\circ\beta=id_{\cp}. \ \
\Box$
%\end{proof}

%\vspace{-0.2cm}
\begin{prop}\label{asso.}Let $D$ be a  regular connection on
$\pi^{-1}(TM)$ with horizontal map$\, \beta$.
\begin{description}
    \item[(a)] The vector field $S$ on $TM$ defined by $S=\beta\circ
    \overline{\eta}\,$ is a semispray.
    \item[(b)]The vector $1$-form $\Gamma$ on $TM$ defined by \vspace{-0.2cm}
    $$\Gamma=2\beta\circ\rho-I\vspace{-0.25cm}$$
    is  a nonlinear connection on $M$.
\vspace{-0.25cm}
\end{description}
This nonlinear connection is characterized by the fact that it
    has the same horizontal and vertical projectors as $D$:
    $h_{\Gamma}=h_{D}=\beta\circ \rho,\, v_{\Gamma}=v_{D}=I-\beta\circ
    \rho$.
\end{prop}

%\begin{proof}
\prof The proof is clear and we omit it. \ \  $\Box$
%\end{proof}

\begin{defn}Let $D$ be a regular connection on
$\pi^{-1}(TM)$ with horizontal map $\beta$.\\
\textbf{--} The semispray  $S=\beta\circ\overline{\eta}$ will be called the semispray associated with $D$.\\
\textbf{--} The nonlinear connection $\Gamma=2\beta\circ\rho-I$ will
be called the nonlinear connection associated with $D$.
\end{defn}

\begin{rem}\label{re} Let $D$ be a  regular connection on
$\pi^{-1}(TM)$ whose horizontal map is $\beta$. The semispray $S$
associated with $D$ coincides with the semispray associated with
$\Gamma$ in the sense of Grifone {\em{\cite{r21}}}. In fact,
  $hS'=(\beta\circ\rho)S'=\beta(\rho S')=\beta\overline{\eta}=S$, where  $S'$ is any arbitrary  semispray.
\end{rem}

\begin{prop}\label{eqv.}Let $(M,L)$ be a Finsler manifold. Let  ${D}$ be a regular connection
on $\pi^{-1}(TM)$ whose connection map is $K$ and whose horizontal
map is $\beta$. Then, the following assertions are equivalent:
\vspace{-0.2cm}
 \begin{description}
    \item[(a)] The (h)hv-torsion  ${T}$ of ${D}$ has the property that
    ${T}( \overline{X},\overline{\eta})=0$,

    \item[(b)] $K=\gamma^{-1}$ on  $V(TM)$,

    \item[(c)] $\widetilde{\Gamma} :=\beta\circ\rho - \gamma\circ K$ is a nonlinear
    connection on $M$.\vspace{-0.2cm}
 \end{description}
 \par Consequently,  if any one of the above assertions
    holds, then $\widetilde{\Gamma}$
    coincides with the nonlinear connection associated with $D$: $\widetilde{\Gamma}=\Gamma=
    2\beta\circ\rho-I$, and in  this case $h_{\Gamma}=h_{D}=\beta\circ\rho$
    and  $\,v_{\Gamma}=v_{D}=\gamma\circ K$.
\end{prop}

%\begin{proof}
\prof\\
\textbf{(a) $\Longleftrightarrow$(b)}: As
$\,S={\beta}\overline{\eta}\,\,$ is a semispray on $M$,
$\,\rho\circ\beta=id_{\cp}\,\, and \,\,\rho\circ\gamma=0\,$, we
have, for all $\overline{X}\in\cp$,\vspace{-0.2cm}
 $$T(\overline{X},\overline{\eta})=
   {D}_{\gamma \overline{X}}\rho({\beta}\overline{\eta})-
  {D}_{{\beta}\overline{\eta}}\rho(\gamma \overline{X})
  - \rho[\gamma \overline{X}, {\beta}\overline{\eta}]
   = {D}_{\gamma \overline{X}}\overline{\eta}
  - \rho[\gamma \overline{X}, S],
  =({K}\circ\gamma) \overline{X}- \overline{X},\vspace{-6pt}$$
  by Lemma\, \ref{iso}.
Consequently,
\begin{equation}\label{eq.n1}
  {T}( \overline{X}, \overline{\eta})= ({K}\circ\gamma- id_{\p}) \overline{X}.
\end{equation}
From which
\begin{equation}\label{eq.n2}
\gamma {T}( \overline{X}, \overline{\eta})=( \gamma\circ K- I)\gamma
 \overline{X}.
\end{equation}
 The result follows from (\ref{eq.n1}) and (\ref{eq.n2}).\vspace{7pt}\\
\textbf{(b) $\Longleftrightarrow$(c)}: If $K=\gamma^{-1}$ on $V(\T
M)$,
 then\vspace{-0.2cm}
   $$J\widetilde{\Gamma} =(\gamma\circ
\rho)\circ(\beta\circ\rho - \gamma\circ K)=\gamma\circ(\rho
\circ\beta)\circ\rho - \gamma\circ(\rho\circ\gamma)\circ K= \gamma
\circ\rho=J,\vspace{-0.2cm}$$
$$\widetilde{\Gamma} J =(\beta\circ\rho - \gamma\circ
K)\circ(\gamma\circ\rho)=\beta\circ (\rho\circ\gamma )\circ\rho -
\gamma\circ (K\circ\gamma )\circ\rho\,= -\gamma
\circ\rho=-J.\vspace{-2pt}$$
 Hence, $\widetilde{\Gamma}$ is a nonlinear
connection.
\par
Conversely, if $\widetilde{\Gamma}$ is a nonlinear connection, we
have from the last relation \vspace{-0.15cm}
$$ (\gamma\circ K )\circ J=J.\vspace{-0.25cm}$$ Hence,
$\gamma\circ K =id_{V(\T M)}$.\\
Similarly,
$$ \gamma \circ \rho =(\gamma\circ K
)\circ(\gamma\circ\rho)=\gamma\circ(\,K \circ\gamma)\circ\rho.$$
From which, since $\gamma :\pi^{-1}(T M)  \To V(\T M)$ is an
isomorphism, $K\circ\gamma =id_{\p}$.
\par
\smallskip
If any one of the assertions \textbf{(a)}-\textbf{(c)} holds,
then\vspace{-0.2cm}
   \begin{align}\label{eq.n3}
    K\circ\gamma =id_{\p}, \quad  \quad
       \gamma\circ K =\left\{
                                \begin{array}{ll}
                                          0 ,& {\,\, on\,\,  H(\T M)} \\
                                         id_{V(\T M)} ,& {\,\, on\,\,    V(\T M)}
                                       \end{array}
                                     \right.\vspace{-0.2cm}
\end{align}
From (\ref{fh1}) and (\ref{eq.n3}), we conclude that\vspace{-0.2cm}
\begin{equation}\label{eq.n4}
\beta\circ\rho +\gamma\circ K =I.\vspace{-0.15cm}
\end{equation}
Consequently, $\widetilde{\Gamma} =\beta\circ\rho - \gamma\circ
K=\beta\circ\rho - (I-\beta\circ\rho)=2\beta\circ\rho - I=\Gamma$,
 which completes the proof. \ \ $\Box$

%\end{proof}

\vspace{7pt}
 We conclude this section by the following lemma which will be used in the sequel.\vspace{-0.2cm}
\begin{lem}\label{bracket} Let $D$ be a regular connection on $\p$
whose (h)hv-torsion tensor $T$ has the property that
$\,T(\overline{X}, \overline{\eta})=0$. Then, we
have{\em:}\vspace{-0.2cm}
   \begin{description}
 \item[(a)] $[\beta \overline{X},\beta \overline{Y}]=
     \gamma\widehat{R}(\overline{X},\overline{Y})
     + \beta(D_{\beta \overline{X}}\overline{Y}-
     D_{\beta \overline{Y}}\overline{X}-Q(\overline{X},\overline{Y})),$

    \item[(b)] $[\gamma \overline{X},\beta \overline{Y}]=-
     \gamma(\widehat{P}(\overline{Y},\overline{X})+D_
     {\beta \overline{Y}}\overline{X})
     +\beta( D_{\gamma \overline{X}}\overline{Y}-T(\overline{X},\overline{Y})),$

   \item[(c)] $[\gamma \overline{X},\gamma \overline{Y}]=
     \gamma(D_{\gamma \overline{X}}\overline{Y}-
     D_{\gamma \overline{Y}}\overline{X}+\widehat{S}(\overline{X},\overline{Y}))$.

     \end{description}
\end{lem}

%\begin{proof}
\prof It should first be noted that, as $D$ is regular and
$T(\overline{X},\overline{\eta})=0$, we have $h=\beta\circ\rho$,
$\,v=\gamma\circ K$, $K\circ\gamma=id_{\cp} $ (cf. Proposition
\ref{eqv.}). \\
We prove only the first part; the other parts can be proved
similarly.
\begin{eqnarray*}
  [\beta \overline{X},\beta \overline{Y}] &=& \gamma
(K
  [\beta \overline{X},\beta \overline{Y}])
  +\beta(\rho [\beta \overline{X},\beta \overline{Y}]) \\
  &=& \gamma
(D_{
  [\beta \overline{X},\beta \overline{Y}]}\overline{\eta})+
  \beta(\rho [\beta \overline{X},\beta \overline{Y}]) \\
  &=&
     \gamma(\widehat{R}(\overline{X},\overline{Y})-D_
     {\beta \overline{Y}}D_
     {\beta\overline{X}}\overline{\eta}+D_
     {\beta\overline{X}}D_
     {\beta \overline{Y}}\overline{\eta})+\beta( D_{\beta
\overline{X}}\overline{Y}-D_{\beta
\overline{Y}}\overline{X}-Q(\overline{X},\overline{Y}))\\
&=&  \gamma\widehat{R}(\overline{X},\overline{Y})
     + \beta(D_{\beta \overline{X}}\overline{Y}-
     D_{\beta \overline{Y}}\overline{X}-Q(\overline{X},\overline{Y})).\ \ \Box
\end{eqnarray*}
%\end{proof}

%%%%%%%%%%%%%%%%%%%%%%%%%%%%%%%%%%%%%%%%%%%%%% SECTION 3. Cartan Connection %%%%%%%%%%%%%%%%%%%%%%%%%%%%%%%%%%%%%%%%%%%%%%%%%%%%%%%%%%%%%%%%

\Section{Cartan Connection}

\par
The aim of the present section is to provide an intrinsic proof of
an intrinsic version of the existence and uniqueness theorem for the
Cartan connection. Moreover, the spray and nonlinear connection
associated with the Cartan connection are investigated.

\begin{thm}  \label{c.ba.}Let $(M,L)$ be a Finsler manifold and
 $g$ the Finsler metric defined by$\,L\,.$ Let $\nabla$ be a regular connection
 on $\pi^{-1}(TM)$ such that
\begin{description}
  \item[(a)]  $\nabla$ is  metric\,{\em:} $\nabla g=0$,

  \item[(b)] The horizontal torsion of $\nabla$ vanishes\,{\em:} $Q=0
  $,
  \item[(c)] The mixed torsion $T$ of $\nabla$ satisfies \,
  $g(T(\overline{X},\overline{Y}), \overline{Z})=g(T(\overline{X},\overline{Z}),\overline{Y})$.\vspace{-0.2cm}
\end{description}
 \par
 Then, the nonlinear connection ${\Gamma}$ associated with
$\nabla$ coincides with the Barthel connection {\em:}
$\Gamma=[J,G]$.
\end{thm}

To prove this theorem, we need the following three lemmas\,:

\vspace{-0.2cm}

\begin{lem}\label{lem.3}Let $(M,L)$ be a Finsler manifold and  $g$ the Finsler metric
defined by$\,L$. For every linear connection $D$ in $\pi^{-1}(TM)$
with torsion  tensor $\textbf{T}$ and curvature tensor $\textbf{K}$,
we have\,:\vspace{-0.1cm}
\begin{description}
    \item[(a)]$\mathfrak{S}_{X,Y,Z}\{\textbf{K}(X,Y)\rho Z+D_{X}\textbf{T}(Y,Z)
    +\textbf{T}(X,[Y,Z])\}=0$.\vspace{-0.1cm}
   \end{description}
If, moreover, $D$ is  metric, then\vspace{-0.1cm}
\begin{description}
    \item[(b)]$g(\textbf{K}(X,Y)\overline{Z},\overline{W})+g(\textbf{K}(X,Y)\overline{W},\overline{Z})=0$.
\end{description}
\end{lem}

\begin{lem} \label{lem.1} Let $(M,L)$ be a Finsler manifold.
Let $D$ be a regular connection on $\pi^{-1}(TM)$ such that
\begin{description}
  \item[(a)]  ${D}$ is vertically  metric\,{\em:} ${D}_{\gamma \overline{X}} g=0$,

  \item[(b)]  The (h)hv-torsion tensor $T$ of ${D}$ satisfies \,
  $g({T}(\overline{X},\overline{Y}), \overline{Z})=
  g({T}(\overline{X},\overline{Z}),\overline{Y})$.
 \end{description}
Then, the (h)hv-torsion tensor $T$ has the property that
$T(\overline{X},\overline{\eta})=0$.
\end{lem}

%\begin{proof}
 \prof If ${D}$ is
a non-metric linear connection on $\pi^{-1}(T M)$ with
 nonzero torsion ${\textbf{T}}$, one can show that
 ${D}$  is completely  determined  by the relation
\vspace{-0.2cm}
\begin{equation}\label{c1eq.r}
     \left.
    \begin{array}{rcl}
    2g({D} _{X}\rho Y,\rho Z)& =& X\cdot g(\rho Y,\rho
                Z)+ Y\cdot g(\rho Z,\rho X)-Z\cdot g(\rho X,\rho
                Y) \\
        & &-g(\rho X,{\textbf{T}}(Y,Z))+g(\rho Y,{\textbf{T}}(Z,X))
        +g(\rho Z,{\textbf{T}}(X,Y)) \\
        & &-g(\rho X,\rho [Y,Z])+g(\rho Y,\rho [Z,X])+g(\rho Z,\rho [X,Y])\\
        & &- ({D} _{X} g)(\rho Y,\rho
                Z)-({D} _{Y} g)(\rho Z,\rho X)+({D} _{Z}g)(\rho X,\rho
                Y).\vspace{-0.2cm}
   \end{array}
  \right\}
 \end{equation}
 for all $X,Y,Z\in\cpp$.
The connection $D$ being regular, let  $h$ and $v$ be the horizontal
and vertical projectors associated with the
 decomposition (\ref{proj.}): $h=\beta\circ\rho$, $ v=I-\beta\circ\rho $.
\par
 Replacing $X, Y, Z$ by $\gamma \overline{X}, {h}Y, {h}Z $
 in (\ref{c1eq.r})
 and using hypotheses  (\textbf{a)}, \textbf{(b)},  taking into account
  the fact that $\rho \circ \gamma=0$ and $\rho\circ h=\rho$, we get\vspace{-0.2cm}
\begin{equation}\label{eq}
2g({D} _{\gamma \overline{X}}\rho Y,\rho Z) =\gamma
\overline{X}\cdot g(\rho Y,\rho Z)+
     g(\rho Y,\rho [hZ,\gamma \overline{X}])+g(\rho Z,\rho [\gamma \overline{X},hY]).\vspace{-0.3cm}
\end{equation}
Now,  \vspace{-0.2cm}
\begin{eqnarray*}
  2\,g(T(\overline{X},\overline{\eta}),\overline{Z})&=&
  2\,g(\textbf{T}(\gamma\overline{X},\beta\overline{\eta}),\overline{Z})\\
     &=&2\,g(D_{\gamma\overline{X}}\overline{\eta},\overline{Z})
  -2\,g(\rho[\gamma\overline{X}, \beta\overline{\eta}],\overline{Z}) \vspace{-0.2cm}.
\end{eqnarray*}
Then, from (\ref{eq}), we get\vspace{-0.2cm}
$$2\,g(T(\overline{X},\overline{\eta}),\overline{Z})=
\gamma\overline{X}\cdot g( \overline{\eta}, \overline{Z})+ g(
\overline{\eta},\rho [\beta\overline{Z},\gamma\overline{X}])-g(
\overline{Z},\rho [\gamma\overline{X},\beta\overline{\eta}]) .$$
Using Lemma \ref{iso}, taking into account the fact that $\beta
\overline{\eta}$ is a semispray, we obtain\vspace{-0.2cm}
$$2\,g(T(\overline{X},\overline{\eta}),\overline{Z})=
\gamma\overline{X}\cdot g( \overline{\eta}, \overline{Z})+ g(
\overline{\eta},\rho [\beta\overline{Z},\gamma\overline{X}])-g(
\overline{Z},\overline{X}) \vspace{-0.2cm}.$$

 Finally, one can show that the sum of the first
two terms on the right-hand side is equal to $g(
\overline{X},\overline{Z})$, from which the result. \ \ $\Box$
%\end{proof}

\begin{lem}\label{lem.4}Under the hypotheses of theorem \ref{c.ba.}, we have\,: \vspace{-0.2cm}
\begin{description}
    \item[(a)]The (h)hv-torsion $T$ is symmetric and $T(\overline{X},\overline{\eta})=0$,

 \item[(b)] The (v)hv-torsion $\widehat{P}$
 is symmetric and
 $\widehat{P}(\overline{X},\overline{\eta})=0$,
  \item[(c)]  The (v)v-torsion $\widehat{S}$ vanishes.
\end{description}
\end{lem}

%\begin{proof}
\prof \\
 \textbf{(a)} By Lemma \ref{lem.1}, taking hypotheses \textbf{(a)} and \textbf{(c)} of Theorem \ref{c.ba.} into account, we have
   $T(\overline{X},\overline{\eta})=0$. It remains to show that $T$ is
   symmetric.\\
    Firstly, one can easily show that
          \begin{equation}\label{eq1}
           g((\nabla_{W}T)(\overline{X},\overline{Y}),\overline{Z})
    =g((\nabla_{W}T)(\overline{X},\overline{Z}),\overline{Y}).
          \end{equation}
Using Lemma \ref{lem.3}\textbf{(a)} for
   $X=\gamma \overline{X},\  Y=\gamma \overline{Y}$ and
$Z=\beta \overline{Z}$, we get\vspace{-0.1cm}
\begin{eqnarray*}
   S(\overline{X}, \overline{Y})\overline{Z}&=&\nabla_{\gamma \overline{Y}}T(\overline{X},\overline{Z} )
   -\nabla_{\gamma \overline{X}}T(\overline{Y},\overline{Z} )
   -\nabla_{\beta \overline{Z}}\textbf{T}(\gamma \overline{X},\gamma \overline{Y} ) -\\
   &&-\textbf{T}(\gamma \overline{X},[\gamma \overline{Y},\beta \overline{Z}]) +\textbf{T}(\gamma \overline{Y},[\gamma \overline{X},\beta \overline{Z}])
   +\textbf{T}([\gamma \overline{X},\gamma \overline{Y}],\beta \overline{Z}). \vspace{-0.1cm}
\end{eqnarray*}
 Now, from  Lemma \ref{bracket} and the fact
 that $\textbf{T}(\gamma \overline{X},\gamma \overline{Y} )=Q(\overline{X},\overline{Y} )=0$, the above equation  reduces to
\begin{equation}\label{La}
\left.
    \begin{array}{rcl}
S(\overline{X},\overline{Y})\overline{Z} &=&
 (\nabla_{\gamma \overline{Y}}T)(\overline{X},\overline{Z})-
 (\nabla_{\gamma \overline{X}}T)(\overline{Y},\overline{Z})
 +T( \widehat{S}(\overline{X}, \overline{Y}),\overline{Z})\\
&
&+T(\overline{X},T(\overline{Y},\overline{Z}))-T(\overline{Y},T(\overline{X},\overline{Z})).\vspace{-0.2cm}
    \end{array}
  \right.
\end{equation}
From which, since
$g(T(\overline{X},\overline{Y}),\overline{Z})=g(T(\overline{X},\overline{Z}),\overline{Y})$,
we have\vspace{-0.2cm}
\begin{equation*}\label{3.eq.7}
\left.
    \begin{array}{rcl}
S(\overline{X},\overline{Y},\overline{Z},\overline{W}) &=&
g((\nabla_{\gamma \overline{Y}}T)(\overline{X},\overline{Z}),
\overline{W}) -g((\nabla_{\gamma
\overline{X}}T)(\overline{Y},\overline{Z}),\overline{W})+
\\
& &+g(T(\overline{X},\overline{W}),T(\overline{Y},\overline{Z}))-
g(T(\overline{Y},\overline{W}),T(\overline{X},\overline{Z}))+\\
&&+g(T( \widehat{S}(\overline{X},
\overline{Y}),\overline{Z}),\overline{W}).\vspace{-0.2cm}
    \end{array}
  \right.
\end{equation*}
Similarly,
\begin{equation*}\label{3.eq.7}
\left.
    \begin{array}{rcl}
S(\overline{X},\overline{Y},\overline{W},\overline{Z}) &=&
g((\nabla_{\gamma \overline{Y}}T)(\overline{X},\overline{W}),
\overline{Z}) -g((\nabla_{\gamma
\overline{X}}T)(\overline{Y},\overline{W}),\overline{Z})+
\\
& &+g(T(\overline{X},\overline{Z}),T(\overline{Y},\overline{W}))
-g(T(\overline{Y},\overline{Z}),T(\overline{X},\overline{W}))+\\
&&+g(T( \widehat{S}(\overline{X},
\overline{Y}),\overline{W}),\overline{Z}).\vspace{-0.2cm}
    \end{array}
  \right.
\end{equation*}
On the other hand, using Lemma \ref{lem.3}\textbf{(b)}, we get
$$S(\overline{X},\overline{Y},\overline{Z},\overline{W})
=-S(\overline{X},\overline{Y},\overline{W},\overline{Z}).$$
 Hence, the above three
equations, together with (\ref{eq1}), yield\vspace{-0.1cm}
\begin{equation}\label{3.eq.8}
    (\nabla_{\gamma \overline{X}}T)(\overline{Y},\overline{Z})-(\nabla_{\gamma
\overline{Y}}T)(\overline{X},\overline{Z})=T(
\widehat{S}(\overline{X},
\overline{Y}),\overline{Z}).\vspace{-0.1cm}
\end{equation}
 Now, setting $\overline{Z}=\overline{\eta}$ in (\ref{3.eq.8}),
 noting that $T(\overline{X},\overline{\eta})=0$, we deduce that
 $T(\overline{X},\overline{Y})=T(\overline{Y},\overline{X})$.

\vspace{4pt}
 \noindent \textbf{(b)}
Using Lemma \ref{lem.3}\textbf{(a)} for $X=\beta \overline{X},\
Y=\gamma \overline{Y}$ and $Z=\beta \overline{Z}$, taking into
account Lemma \ref{bracket} and the fact that $Q=0$, we get
\begin{equation*}
\left.
    \begin{array}{rcl}
P(\overline{X},\overline{Y})\overline{Z}-P(\overline{Z},\overline{Y})\overline{X}&=&
(\nabla_{\beta
\overline{Z}}T)(\overline{Y},\overline{X})-(\nabla_{\beta
\overline{X}}T)(\overline{Y},\overline{Z})-\\&
&-T(\widehat{P}(\overline{Z},\overline{Y}),\overline{X})+
T(\widehat{P}(\overline{X},\overline{Y}),\overline{Z}).
    \end{array}
  \right.\vspace{-0.2cm}
\end{equation*}
From which, making use of hypothesis \textbf{(c)} of Theorem
\ref{c.ba.} and  the symmetry of $T$, we have\vspace{-0.2cm}
\begin{equation*}
\left.
    \begin{array}{rcl}
P(\overline{X},\overline{Y},\overline{Z},\overline{W})-
P(\overline{Z},\overline{Y},\overline{X}, \overline{W})&=&
g((\nabla_{\beta\overline{Z}}T)(\overline{Y},\overline{X}),
\overline{W})-g((\nabla_{\beta
\overline{X}}T)(\overline{Y},\overline{Z}), \overline{W})\\
&-&g(T(\overline{X},\overline{W}),\widehat{P}(\overline{Z},\overline{Y}))+
g(T(\overline{Z},\overline{W}),\widehat{P}(\overline{X},\overline{Y})).
    \end{array}
  \right.\vspace{-0.2cm}
\end{equation*}
By cyclic permutation on $\overline{X},\overline{Z},\overline{W}$ of
the above equation, on gets
\begin{equation*}
\left.
    \begin{array}{rcl}
P(\overline{W},\overline{Y},\overline{X},\overline{Z})-
P(\overline{X},\overline{Y},\overline{W}, \overline{Z})&=&
g((\nabla_{\beta\overline{X}}T)(\overline{Y},\overline{W}),
\overline{Z})-g((\nabla_{\beta
\overline{W}}T)(\overline{Y},\overline{X}), \overline{Z})\\
&-&g(T(\overline{W},\overline{Z}),\widehat{P}(\overline{X},\overline{Y}))+
g(T(\overline{X},\overline{Z}),\widehat{P}(\overline{W},\overline{Y}))
    \end{array}
  \right.\vspace{-0.3cm}
\end{equation*}
and
\begin{equation*}
\left.
    \begin{array}{rcl}
P(\overline{Z},\overline{Y},\overline{W},\overline{X})-
P(\overline{W},\overline{Y},\overline{Z}, \overline{X})&=&
g((\nabla_{\beta\overline{W}}T)(\overline{Y},\overline{Z}),
\overline{X})-g((\nabla_{\beta
\overline{Z}}T)(\overline{Y},\overline{W}), \overline{X})\\
&-&g(T(\overline{Z},\overline{X}),\widehat{P}(\overline{W},\overline{Y}))+
g(T(\overline{W},\overline{X}),\widehat{P}(\overline{Z},\overline{Y})).
    \end{array}
  \right.\vspace{-0.3cm}
\end{equation*}
Making use of the above three relations, together with the identity
$P(\overline{X},\overline{Y},\overline{Z},\overline{W})
=-P(\overline{X},\overline{Y}, \overline{W},\overline{Z})$ (by Lemma
\ref{lem.3}\textbf{(b)}), it follows that \vspace{-0.1cm}
 \begin{equation}\label{01}
\left.
    \begin{array}{rcl}
   P(\overline{X},\overline{Y},\overline{Z},\overline{W})&=&
   g((\nabla_{\beta\overline{Z}}T)(\overline{Y},\overline{X}),
\overline{W})
   -g((\nabla_{\beta
\overline{W}}T)(\overline{Y},\overline{X}), \overline{Z})
   \\
    &&-g(T(\overline{X},\overline{W}),\widehat{P}(\overline{Z},\overline{Y}))
   +g(T(\overline{X},\overline{Z}),\widehat{P}(\overline{W},\overline{Y})).
\end{array}
  \right.
\end{equation}
\par Setting $
    \overline{X}=\overline{Z}=\overline{\eta}$ in (\ref{01}), taking into account the fact  that
    $T(\overline{X}, \overline{\eta})=0$, ${K}\circ {\beta}\,=0$
     and that the metric tensor $g$ is
    nongegenerate, we obtain $\widehat{{P}}(\overline{\eta}, \overline{Y})=0$ for all $\overline{Y} \in \cp$.
   Consequently, Equation (\ref{01}) for $\overline{Z}=\overline{\eta}$
   implies that
\begin{equation}\label{02}
   \widehat{P}(\overline{X}, \overline{Y})=
   (\nabla_{\beta\overline{\eta}}T)(\overline{X},\overline{Y}),
\end{equation}
where we have used (\ref{eq1}). The symmetry of $\widehat{P}$
follows then from the  symmetry of $T$.

\vspace{4pt}
 \noindent \textbf{(c)} From  (\ref{La}) and (\ref{3.eq.8}), we get
\begin{equation*}
S(\overline{X},\overline{Y})\overline{Z}=
T(\overline{X},T(\overline{Y},\overline{Z}))-T(\overline{Y},T(\overline{X},\overline{Z})).\vspace{-0.2cm}
\end{equation*}
Setting $\overline{Z}=\overline{\eta}$ and noting that
$T(\overline{X},
 \overline{\eta})=0$, the result follows. \ \ $\Box$
%\end{proof}

\vspace{12pt}

%\begin{proof}
\noindent\textit{\textbf{Proof of Theorem \ref{c.ba.}}}\,:
\par
 As $ T(\overline{X}, \overline{\eta})=0$ for the
 connection $\nabla $ (by Lemma \ref{lem.4}), then  by Proposition \ref{eqv.}, it follows
$K=\gamma^{-1}$ on $V(\T M)$ and the associated  nonlinear
connection $\Gamma$ is given by $\Gamma=\beta\circ\rho -\gamma\circ
K$.
\par We prove that $\Gamma$ enjoys the following properties:

\vspace{4pt}
 \noindent\textit{\textbf{$\Gamma $ is conservative $(d_{h}E=0)$:}}
\vspace{-0.3cm}
       $$d_{h} E(X)= i_{h}dE(X)=hX\cdot E
        =\frac{1}{2} hX \cdot g(\overline{\eta}, \overline{\eta}),
        =g( \nabla_{hX}\overline{\eta},\overline{\eta})= 0.\vspace{-0.3cm}$$

\vspace{4pt}
 \noindent \textit{\textbf{$\Gamma$ is
homogenous $(\,[\mathcal{C},\Gamma]=0)$:}}
\par
It is easy to show that
  $$
    [\mathcal{C},v]X=-v[\mathcal{C},hX].
  $$
As $v=\gamma\circ K$, $h=\beta\circ\rho$ and
$\gamma\circ\overline{\eta}=\mathcal{C}$, then \vspace{-3pt}
$$[\mathcal{C},v]X=-(\gamma \circ K)[\gamma \overline{\eta},\beta\rho X].\vspace{-3pt}$$
 Now, by Lemma \ref{bracket}\textbf{(b)}, noting that
$\widehat{P}( \overline{X},\overline{\eta})=0$ (Lemma \ref{lem.4}),
$K\circ\gamma=id_{\cp}$ and $ K\circ \beta=0$, we obtain
\begin{eqnarray*}
    [\mathcal{C},v]X
    &=&-(\gamma\circ  K)\{-
     \gamma(\widehat{P}(\rho X,\overline{\eta})+\nabla_
     {\beta \rho X}\overline{\eta})
     +\beta( \nabla_{\gamma \overline{\eta}}\rho X-T(\overline{\eta},\rho
     X))\}\\
     &=&\gamma \{
     (\widehat{P}(\rho X,\overline{\eta})+\nabla_
     {\beta \rho X}\overline{\eta})\}=0.
  \end{eqnarray*}
Consequently, $ [\mathcal{C},\Gamma]= -2[\mathcal{C},v]=0$.

\vspace{4pt}
 \noindent
\textit{\textbf{$\Gamma$ is torsion-free $(\,[J,\Gamma]=0)$:}}
\begin{eqnarray*}
  [J,v](X,Y) &=& [JX,vY]+ [vX,JY]+vJ[X,Y]+Jv[X,Y]\\
   && -J[vX,Y]-J[X,vY]-v[JX,Y]-v[X,JY].
\end{eqnarray*}
As $J\circ v=0$, $v\circ J=J$ and the vertical distribution is
completely integrable, we get
\begin{eqnarray*}
[J,v](X,Y) &=&J[hX,hY]-v[JX,hY]-v[hX,JY]\\
   &=&J[\beta \rho X,\beta \rho Y]-v[\gamma \rho X,\beta \rho Y]+v[\gamma \rho Y,\beta \rho
   X].
\end{eqnarray*}
From which, together with Lemma \ref{bracket},  taking into account
the fact that $Q=0$, we obtain
  \begin{eqnarray*}
  [J,v](X,Y)&=&J\{\gamma\widehat{R}(\rho{X},\rho{Y})
     + \beta(\nabla_{h{X}}\rho{Y}-
     \nabla_{h{Y}}\rho{X})\}\\
    &&-(\gamma\circ K)\{-
     \gamma(\widehat{P}(\rho{Y},\rho{X})+\nabla_
     {h{Y}}\rho{X})
     +\beta( \nabla_{J{X}}\rho{Y}-T(\rho{X},\rho{Y}))\}\\
     &&+(\gamma\circ K)\{-
     \gamma(\widehat{P}(\rho{X},\rho{Y})+\nabla_
     {h{X}}\rho{Y})
     +\beta( \nabla_{J{Y}}\rho{X}-T(\rho{Y},\rho{X}))\}.
     \end{eqnarray*}
Noting that $\,J\circ\gamma=0$, $\,J\circ\beta=\gamma$,
$\,K\circ\gamma=id_{\cp}$, $\,K\,\circ\beta=0$
        and that
    $ \widehat{P}$ is symmetric (by Lemma \ref{lem.4}(c)), it follows that $[J,v]=0$.
    From which
    $t:=\frac{1}{2}[J,\Gamma]=-[J,v]=0$.
   \par From the above consideration, $\Gamma=\beta\circ\rho -\gamma\circ K$
   is a conservative torsion-free homogenous
   nonlinear connection. By the uniqueness of the Barthel connection
    (Theorem \ref{th.9a}), it follows that  $\Gamma $ coincides with the Barthel
   connection $[J,G]$.\ \ $\Box$
%\end{proof}

\vspace{4pt} In view Theorem \ref{c.ba.} and Remark \ref{re}, we
have the\vspace{-0.3cm}
\begin{cor}\label{Cth.1}
 The semispray associated with  the connection $\nabla$  (of Theorem \ref{c.ba.})
 is a spray which coincides with the
canonical spray.
\end{cor}

\begin{rem} From Theorem \ref{c.ba.}, Proposition \ref{asso.}, and Equation (\ref{eq.n4}),
the nonlinear connection associated with the connection $\nabla$ can
be expressed in different equivalent forms:
\begin{equation}\label{eq.n5}
  \Gamma=2\,\beta\circ \rho-I=I-2\,\gamma \circ K=\beta\circ \rho -\gamma\circ
  K=[J,G],\vspace{-0.2cm}
\end{equation}
which provides a strong link between the KG-approach and the
PB-approach.
\end{rem}

 \par Now, we have the following fundamental
result\,:\vspace{-0.2cm}
\begin{thm} \label{th.1}Let $(M,L)$ be a Finsler manifold and  $g$  the Finsler metric
defined by $L$. There exists a unique regular connection $\nabla$ on
$\pi^{-1}(TM)$ such that\vspace{-0.2cm}
\begin{description}
  \item[(a)]  $\nabla$ is  metric\,{\em:} $\nabla g=0$,

  \item[(b)] The horizontal torsion of $\nabla$ vanishes\,{\em:} $Q=0
  $,
  \item[(c)] The mixed torsion $T$ of $\nabla$ satisfies \,
  $g(T(\overline{X},\overline{Y}), \overline{Z})=g(T(\overline{X},\overline{Z}),\overline{Y})$.\vspace{-0.2cm}
\end{description}
 \par
 Such a connection is called the Cartan
connection associated with  the Finsler manifold $(M,L)$.

\end{thm}

%\begin{proof}
\prof  The connection $\nabla$ being regular, let $h$ and $v$ be its
 horizontal and vertical projectors. Then, by  Theorem \ref{c.ba.}, these projectors
  coincide with the corresponding projectors of the Barthel connection.
  \par
   First we prove the {\it \textbf{uniqueness}}. As $\nabla$ is
a metric linear  connection on $\pi^{-1}(T M)$ with nonzero torsion
${\textbf{T}}$, then, by (\ref{c1eq.r}),
 $\nabla$  is completely  determined by the relation
\vspace{-0.2cm}
\begin{equation}\label{ceq.r}
     \left.
    \begin{array}{rcl}
    2g(\nabla _{X}\rho Y,\rho Z)& =& X\cdot g(\rho Y,\rho
                Z)+ Y\cdot g(\rho Z,\rho X)-Z\cdot g(\rho X,\rho
                Y) \\
        & &-g(\rho X,{\textbf{T}}(Y,Z))+g(\rho Y,{\textbf{T}}(Z,X))
        +g(\rho Z,{\textbf{T}}(X,Y)) \\
        & &-g(\rho X,\rho [Y,Z])+g(\rho Y,\rho [Z,X])+g(\rho Z,\rho [X,Y])
        ,\vspace{-0.2cm}
   \end{array}
  \right\}
 \end{equation}
 for all $X, Y,Z\in\cpp$.
\par

  Replacing $X, Y, Z$ by $hX, hY, hZ $
in (\ref{ceq.r})
 and using axiom \textbf{(b)} and the fact that $\rho\circ h=\rho$, we get\vspace{-0.2cm}
\begin{equation}\label{ceq.rt5}
    \left.
    \begin{array}{rcl}
     2g(\nabla _{hX}\rho Y,\rho Z) & = & hX\cdot g(\rho Y,\rho
                Z)+ hY\cdot g(\rho Z,\rho X)-hZ\cdot g(\rho X,\rho
                Y) \\
        & -&g(\rho X,\rho [hY,hZ])+g(\rho Y,\rho
[hZ,hX])+g(\rho Z,\rho
    [hX,hY]).\vspace{-0.2cm}
    \end{array}
  \right.
\end{equation}
Similarly, by replacing $X, Y, Z$ by $vX, hY, hZ $ in (\ref{ceq.r}),
where $vX=\gamma \overline{X}$ for some $\overline{X}\in \cp$,
 and using axiom  \textbf{(c)} and the fact that $\rho \circ v=0$, we get  \vspace{-0.2cm}
              \begin{equation}\label{ceq.rt6}
              2g(\nabla _{vX}\rho Y,\rho Z) =vX\cdot g(\rho Y,\rho Z)+
     g(\rho Y,\rho [hZ,vX])+g(\rho Z,\rho [vX,hY]).\vspace{-0.2cm}
              \end{equation}
\noindent Hence, $\nabla_{X}\rho Y$ is uniquely determined by the
right-hand side of Equations (\ref{ceq.rt5}) and (\ref{ceq.rt6}),
where $h$ and $v$ are known a priori. \vspace{7pt}
\par To prove  the {\it \textbf{existence}}, we define $\nabla$ by
the requirement that (\ref{ceq.rt5}) and (\ref{ceq.rt6}) hold for
all $X, Y, Z \in \x$. Now, we have to prove that the connection
$\nabla$ satisfies the conditions of Theorem \ref{th.1}\,:\smallskip\\
\noindent\textit{\textbf{ $\nabla$ satisfies condition (a)}}\,: By
using (\ref{ceq.rt6}), we get\vspace{-0.2cm}
\begin{eqnarray*}
   2(\nabla _{vX}g)(\rho Y,\rho Z)&=& 2 \{vX\,\cdot g(\rho Y,\rho Z)-
   g(\nabla _{vX}\rho Y,\rho Z)-g(\rho Y,\nabla _{vX}\rho Z)\}  \\
   &=& 2\, vX\,\cdot g(\rho Y,\rho Z)-\\
   &&- \{ vX\cdot g(\rho Y,\rho Z)+
     g(\rho Y,\rho [hZ,vX])+g(\rho Z,\rho [vX,hY])\}\\
     &&- \{ vX\cdot g(\rho Z,\rho Y)+
     g(\rho Z,\rho [hY,vX])+g(\rho Y,\rho [vX,hZ])\}=0.
 \end{eqnarray*}
Similarly, using (\ref{ceq.rt5}), one can show that $ (\nabla
_{hX}g)(\rho Y,\rho Z)=0$.\smallskip\\
\textit{\textbf{ $\nabla$ satisfies condition (b)}}\,:
 From the definition of the (h)h-torsion tensor $Q$ of $ \nabla $,
 using (\ref{ceq.rt5}), we have
\begin{eqnarray*}
  2g(\textbf{T}(hX,hY),\rho Z) &=&2g(\nabla _{hX}\rho Y-\nabla _{hY}\rho X-\rho[hX,hY],  \rho Z)  \\
  &=&\{hX\cdot g(\rho Y,\rho
                Z)+ hY\cdot g(\rho Z,\rho X)-hZ\cdot g(\rho X,\rho
                Y) \\
        & & -g(\rho X,\rho [hY,hZ])+g(\rho Y,\rho
[hZ,hX])+g(\rho Z,\rho
    [hX,hY])\}\\
    &&-\{hY\cdot g(\rho X,\rho
                Z)+ hX\cdot g(\rho Z,\rho Y)-hZ\cdot g(\rho Y,\rho
                X) \\
        & & -g(\rho Y,\rho [hX,hZ])+g(\rho X,\rho
[hZ,hY])+g(\rho Z,\rho
    [hY,hX])\}\\
    &&-2g(\rho[hX,hY],\rho Z).
\end{eqnarray*}
From which, it  follows that $g(\textbf{T}(hX,hY),\rho Z)=0$, for all $X,Y,Z\in \cpp$.\smallskip\\
\textit{\textbf{ $\nabla $ satisfies condition (c)}}\,: By using
(\ref{ceq.rt6}) and the fact that $\nabla$ is metric, we
get\vspace{-0.2cm}
\begin{eqnarray*}
   g(\textbf{T}(vX, hY),\rho Z)&=& g(\nabla _{vX}\rho Y-\rho[vX,hY],\rho Z)\\
   &=& vX\,\cdot g(\rho Y,\rho Z)-g(\rho Y,\nabla _{vX}\rho Z)-g(\rho[vX,hY],\rho Z) \\
    &=& \{vX\,\cdot g(\rho Y,\rho Z)+g(\rho[hY,vX],\rho Z)\}-g(\rho Y,\nabla _{vX}\rho Z) \\
   &=& 2g(\rho Y,\nabla _{vX}\rho Z)-g(\rho[vX,hZ],\rho Y)
    -g(\rho Y,\nabla _{vX}\rho Z)\\
     &=& g(\textbf{T}(vX, hZ),\rho Y).
  \end{eqnarray*}
This completes the proof of Theorem \ref{th.1}.$\ \ \Box$
%\end{proof}

\vspace{7pt}
In view of the above theorem, we have\vspace{-0.2cm}
\begin{thm}\label{le.cc3}The Cartan connection $\nabla$
 is uniquely determined by  the following relations:
   \begin{description}
     \item[(a)]  $ 2g(\nabla _{\gamma\overline{X}}
      \overline{Y}, \overline{Z})
      =\gamma\overline{X}\cdot g( \overline{Y},\overline{Z})+
     g( \overline{Y},\rho [\beta\overline{Z},\gamma\overline{X}])
     +g(\overline{Z},\rho [\gamma\overline{X},\beta \overline{Y}])$,

     \item[(b)]  $ 2g(\nabla _{\beta\overline{X}}\rho Y,\rho Z)
     = \beta\overline{X}\cdot g( \overline{Y},
                \overline{Z})+
                \beta\overline{Y}\cdot g( \overline{Z},\overline{ X})
                -\beta\overline{Z}\cdot g( \overline{X},
                \overline{Y})$\\
  $ { \qquad \qquad \qquad \quad }
  -g( \overline{X},\rho [\beta\overline{Y},\beta\overline{Z}])+g( \overline{Y},\rho
[\beta\overline{Z},\beta\overline{X}])+g( \overline{Z},\rho
    [\beta\overline{X},\beta\overline{Y}])$,
 \end{description}
 where $\beta$ is the horizontal map of the Cartan
 connection {\em(}given by the relation $\beta\circ\rho=h$; $h$ being the horizontal
 projector of Barthel connection{\em)}.
\end{thm}

%%%%%%%%%%%%%%%%%%%%%%%%%%%%%%%%%%%%%%%%%%%%%% SECTION 4: Berwald Connection %%%%%%%%%%%%%%%%%%%%%%%%%%%%%%%%%%%%%%%%%%%%%%%%%%%%%%%%%%%%%%%%%

\Section{Berwald Connection}
 In this section, we provide an intrinsic proof of the
existence and uniqueness theorem for the Berwald connection
$D^{\circ}$. Moreover, we deduce an explicit expression relating
this connection and the Cartan connection $\nabla$.

 \vspace{7pt}
 The following lemma is useful for subsequent use.\vspace{-0.2cm}
\begin{lem} \label{le.n1} Let $D$ be a regular connection on $\pi^{-1}(TM)$ with the following
properties:
\begin{description}
 \item[(a)] $D_{hX}L=0$, $h$ being the horizontal projector of $D$,
 \item[(b)] ${{D}}$ is torsion-free\,{\em:} ${\textbf{T}}=0 $,
 \item[(c)] The (v)hv-torsion tensor $\widehat{{{P}}}$ of ${D}$ vanishes\,:
   $\widehat{{{P}}}(\overline{X},\overline{Y})= 0$.
   \end{description}
Then, the nonlinear connection associated with $D$ coincides with
the Barthel connection.
\end{lem}

%\begin{proof}
\prof The proof is similar to the proof of Theorem \ref{c.ba.},
taking into account Lemma \ref{bracket} together with the given
properties of $D$. \ \ $\Box$
%\end{proof}

\begin{rem}\label{rrem1}
Let $D$ be a regular connection on $\pi^{-1}(TM)$. If the nonlinear
connection associated with $D$ coincides with the Barthel
connection, then the horizontal and vertical projectors of $D$
coincide with the horizontal and vertical projectors of the Cartan
connection $\nabla$ {\em(}Theorem \ref{c.ba.}{\em)}.
\end{rem}

\par Now, we announce the main result of this section, namely, the
existence and uniqueness theorem of the Berwald connection.
\begin{thm} \label{bth2.h2} Let $(M,L)$ be a Finsler manifold. There exists a
unique regular connection ${{D}}^{\circ}$ on $\pi^{-1}(TM)$ such
that
\begin{description}
 \item[(a)] $D^{\circ}_{h^{\circ}X}L=0$,
  \item[(b)]   ${{D}}^{\circ}$ is torsion-free\,{\em:} ${\textbf{T}}^{\circ}=0 $,
  \item[(c)]The (v)hv-torsion tensor $\widehat{{{P}}}^{\circ}$ of ${D}^{\circ}$ vanishes:
   $\widehat{{{P}}}^{\circ}(\overline{X},\overline{Y})= 0$.
  \end{description}
  \par Such a connection is called the Berwald
  connection associated with the Finsler manifold $(M,L)$.
\end{thm}

%\begin{proof}
\prof  First we prove the {\it \textbf{uniqueness}}. As
${{D}}^{\circ}$ is a non-metric linear  connection on $\pi^{-1}(T
M)$ with zero torsion, then, by (\ref{c1eq.r}),
 ${{D}}^{\circ}$  is completely  determined by the relation
\vspace{-0.2cm}
\begin{equation}\label{beq.2}
    \left.
    \begin{array}{rcl}
    2g({{D}}^{\circ}_{X}\rho Y,\rho Z)& =& X\cdot g(\rho Y,\rho
                Z)+ Y\cdot g(\rho Z,\rho X)-Z\cdot g(\rho X,\rho
                Y) \\
        && -g(\rho X,\rho [Y,Z])+g(\rho Y,\rho [Z,X])+g(\rho Z,\rho [X,Y])\\
        && - ({{D}}^{\circ} _{X} g)(\rho Y,\rho
                Z)-({{D}}^{\circ} _{Y} g)(\rho Z,\rho X)+
                ({{D}}^{\circ} _{Z}g)(\rho X,\rho
                Y),\vspace{-0.2cm}
    \end{array}
    \right \}
 \end{equation}
for all $X,Y,Z\in\cpp$.
\par
The connection $D^{\circ}$ being regular, let $h^{\circ}$ and
$v^{\circ}$ be its
 horizontal and vertical projectors.
 According to the axioms of the theorem, the nonlinear connection associated with
  $D^{\circ}$ coincides with the Barthel connection (Lemma \ref{le.n1}).
  Hence, we have $v^{\circ}=v$, $h^{\circ}=h$ ($h$ and $v$ being the horizontal and vertical
  projectors of the Cartan connection (Remmark \ref{rrem1})) and, consequently,
$K^{\circ}=K$, $\beta^{\circ}=\beta$.

\par
By replacing $X, Y, Z$ in (\ref{beq.2}) by $vX, hY, hZ$
respectively, noting that $\rho\circ v=0$ and $\rho\circ h=\rho$, we
get\vspace{-0.3cm}
\begin{equation*}
2g({{D}}^{\circ} _{vX}\rho Y,\rho Z) =vX\cdot g(\rho Y,\rho Z)+
     g(\rho Y,\rho [Z,vX])+g(\rho Z,\rho [vX,Y])-({{D}}^{\circ} _{vX} g)(\rho Y,\rho
                Z).\vspace{-0.3cm}
\end{equation*}
Using  Theorem \ref{le.cc3}\textbf{(a)}, the above equation implies
that \vspace{-0.2cm}
\begin{equation}\label{eq.n6}
 2g({{D}}^{\circ} _{vX}\rho Y,\rho Z)=2g(\nabla _{vX}\rho Y,\rho Z)
 -({{D}}^{\circ} _{vX} g)(\rho Y,\rho
                Z).\vspace{-0.2cm}
\end{equation}
Consequently,
\begin{eqnarray*}
  2g(\textbf{T}^{\circ}({vX}, hY),\rho Z)&=&
  2g({{D}}^{\circ} _{vX}\rho Y-\rho[vX,hY],\rho Z)\\
  &=&2g(\nabla _{vX}\rho Y-\rho[vX,hY],\rho Z)
 -({{D}}^{\circ} _{vX} g)(\rho Y,\rho  Z)\\
   &=& 2g(\textbf{T}({vX}, hY),\rho Z)
 -({{D}}^{\circ} _{vX} g)(\rho Y,\rho  Z).
  \end{eqnarray*}
From which, taking axiom \textbf{(b)} into account, we get   $$
({{D}}^{\circ} _{vX} g)(\rho Y,\rho Z)= 2g(\textbf{T}({vX}, hY),\rho
Z).$$
 Consequently, (\ref{eq.n6}) reduces to
\begin{equation}\label{bqq.02}
 {{D}}^{\circ} _{vX}\rho Y=\nabla _{vX}\rho Y
 -\textbf{T}({vX}, hY).\vspace{-0.1cm}
\end{equation}
\par Similarly, using axiom \textbf{(c)} and  noting that $K\circ J=\gamma$ and
$K\circ h=0$, we get \vspace{-0.2cm}
\begin{eqnarray*}
  0=\widehat{P}^{\circ}(hX,JY)&=&P^{\circ}(hX,JY)\overline{\eta} \\
  &=& -{{D}}^{\circ}_{hX} {D}^{\circ}_{JY}\overline{\eta}
  +{{D}}^{\circ}_{JY} {D}^{\circ}_{hX}\overline{\eta} +{{D}}^{\circ}_{[hX,JY]}\overline{\eta} \\
  &=& - {{D}}^{\circ}_{hX}\rho Y+ K[hX,JY].\vspace{-0.3cm}
  \end{eqnarray*}
From which,
\begin{equation}\label{eq.n7}
    {{D}}^{\circ}_{hX}\rho Y= K[hX,JY].
\end{equation}
Using Lemma \ref{bracket}, (\ref{eq.n7}) may also be written in the
form
  \begin{equation}\label{bqq1.2}
D^{\circ}_{hX}\rho Y=\nabla_
     {hX}\rho Y +\widehat{P}(\rho {X},\rho{Y}).
   \end{equation}
Consequently, from (\ref{bqq.02}) and (\ref{bqq1.2}), the full
expression of $D^{\circ}_{X}\overline{Y}$ is given by
\begin{equation}\label{bqq1.3}
 {{D}}^{\circ}_{X}\overline{Y} = \nabla _{X}\overline{Y}
+ \widehat{P}(\rho X,\overline{Y})-T(K X,
\overline{Y}).\vspace{-0.2cm}
\end{equation}
Hence ${D}^{\circ}_{X}\overline{Y}$ is uniquely determined by the
right-hand side of (\ref{bqq1.3}).

\vspace{7pt}
\par To prove the {\it
\textbf{existence}} of ${D}^{\circ}$, we define ${D}^{\circ}$ by the
requirement that (\ref{bqq1.3}) holds (or, equivalently,
(\ref{bqq.02}) and (\ref{bqq1.2}) hold) for all $X \in \x \ and \
\overline{Y}\in\cp$. Now, we have to prove that the connection
${D}^{\circ}$ satisfies the conditions of Theorem \ref{bth2.h2}:

\vspace{5pt}\noindent \textit{\textbf{ ${D}^{\circ}$ satisfies
condition (a)}}:  Setting $\overline{Y}=\overline{\eta}$ in
(\ref{bqq1.3}),  taking into account the fact that
  $\widehat{P}(\overline{X},\overline{\eta})=0=T(\overline{X},\overline{\eta})$, it follows
that ${K}^{\circ}=K$. As $T^{\circ}=0$, we have
${v}^{\circ}=\gamma\circ K^{\circ}=\gamma\circ K=v$ (by Proposition
\ref{eqv.}). Consequently, $h^{\circ}=h$ and hence \vspace{-0.2cm}
$$0=L\,{D}^{\circ}_{{h^{\circ}}X}L= {D}^{\circ}_{{h^{\circ}}X}E={D}^{\circ}_{{h}X}E=hX\cdot E=d_{h}E(X).$$

\vspace{-1pt}\noindent \textit{\textbf{ ${{D}}^{\circ}$ satisfies
condition (b)}}: From (\ref{bqq1.3}), $\widehat{P}$ being symmetric,
we have \vspace{-0.2cm}
\begin{eqnarray*}
  {\textbf{T}}^{\circ}(X,Y) &=& {{D}}^{\circ}_{X} \rho Y- {{D}}^{\circ}_{Y} \rho X-\rho [X, Y]  \\
   &=&\nabla _{X}\rho {Y}
+ \widehat{P}(\rho X,\rho {Y})-\textbf{T}({vX}, hY)-\\
&&-\nabla _{Y}\rho {X} - \widehat{P}(\rho Y,\rho {X})+\textbf{T}({vY}, hX)-\rho[X,Y]\\
&=&\textbf{T}(X, Y )-\textbf{T}({vX}, hY)+\textbf{T}({vY}, hX)=0.
\end{eqnarray*}

\vspace{-1pt}\noindent\textit{\textbf{ ${D}^{\circ}$ satisfies
condition (c)}}: Using (\ref{bqq1.3}) and the properties of the
Cartan connection  $\nabla$, we get\vspace{-0.2cm}
\begin{eqnarray*}
  \widehat{{P}}^{\circ}(\overline{X},\overline{Y})
  &=&\textbf{{{K}}}^{\circ}(\beta\overline{X},\gamma\overline{Y})\overline{\eta}
  =-{{D}}^{\circ}_{\beta\overline{X}}D^{\circ}_{\gamma\overline{Y}}\overline{\eta}
  +{{D}}^{\circ}_{\gamma\overline{Y}}D^{\circ}_{\beta\overline{X}}\overline{\eta}
+{{D}}^{\circ}_{[\beta\overline{X},\gamma\overline{Y}]}\overline{\eta}\\
&=&-\nabla_{\beta\overline{X}}
\nabla_{\gamma\overline{Y}}\overline{\eta}-
\widehat{P}(\overline{X},\nabla_{\gamma\overline{Y}}\overline{\eta})
+\nabla_{[\beta\overline{X},\gamma\overline{Y}]}\overline{\eta}
\\
&=&\{-\nabla_{\beta\overline{X}}
\nabla_{\gamma\overline{Y}}\overline{\eta}
+\nabla_{[\beta\overline{X},\gamma\overline{Y}]}\overline{\eta}\}-
\widehat{P}(\overline{X},\overline{Y})=0.
\end{eqnarray*}
This complete the proof. \ \ $\Box$
%\end{proof}

\vspace{7pt}
\par
In view of the above theorem, we have\,:
 \vspace{-0.1cm}
\begin{thm}\label{th.5} The Berwald connection $D^{\circ}$ is explicitly expressed in
terms of the Cartan connection $\nabla$ in the form:
 \vspace{-0.1cm}
  \begin{equation}\label{b10}
     {{D}}^{\circ}_{X}\overline{Y} = \nabla _{X}\overline{Y}
+{\widehat{P}}(\rho X,\overline{Y}) -T(K X,\overline{Y}).
\vspace{-0.2cm}
\end{equation}
In particular, we have\vspace{-0.1cm}
\begin{description}
  \item[(a)] $ {{D}}^{\circ}_{\gamma \overline{X}}\overline{Y}=\nabla _{\gamma
  \overline{X}}\overline{Y}-T(\overline{X},\overline{Y})=\rho[\gamma
\overline{X}, \beta \overline{Y}]$.

 \item[(b)] $ {{D}}^{\circ}_{\beta \overline{X}}\overline{Y}=\nabla _{\beta
  \overline{X}}\overline{Y}+\widehat{P}(\overline{X},\overline{Y})=K[\beta
\overline{X}, \gamma \overline{Y}].$
\end{description}
\end{thm}

\begin{rem}
From the above consideration, it should be noted that the semispray
associated with the Berwald connection is a spray which coincides
with the canonical spray. Moreover, the nonlinear connection
associated with the Berwald connection coincides with the Barthel
connection.
\end{rem}

Concerning the metricity properties of $D^{\circ}$, we terminate
with the following result which is not difficult to prove.
\vspace{-0.11cm}
\begin{prop} The Berwald connection $D^{\circ}$ has the properties:
\begin{description}
\item[(a)] $({{D}}^{\circ}_{\gamma
\overline{X}}\,g)(\overline{Y},\overline{Z}) =2 g(T(\overline{X},
\overline{Y}),\overline{Z})$.

\item[(b)] $({{D}}^{\circ}_{\beta \overline{X}}\,g)(\overline{Y},\overline{Z}) =-2
  g(\widehat{P}(\overline{X},\overline{Y}),\overline{Z})$.

\item[(c)] $ {{D}}^{\circ}_{G}\,g=0$.
\end{description}
Consequently,\\
\textbf{--} a Finsler manifold $(M,L)$ is Riemannian if and only if
${{D}}^{\circ}_{\gamma \overline{X}}\,g=0$.\\
\textbf{--} a Finsler manifold $(M,L)$ is Landsbergian if and only
if ${{D}}^{\circ}_{\beta \overline{X}}\,g=0$.
\end{prop}

%%%%%%%%%%%%%%%%%%%%%%%%%%%%%%%%%%%%%%%%%%%%%%%%%%%% REFERENCES %%%%%%%%%%%%%%%%%%%%%%%%%%%%%%%%%%%%%%%%%%%%%%%%%%%%%%%%%%%%%%%%%%%%%%%%%%%

\providecommand{\bysame}{\leavevmode\hbox
to3em{\hrulefill}\thinspace}
\providecommand{\MR}{\relax\ifhmode\unskip\space\fi MR }
\providecommand{\MRhref}[2]{%
  \href{http://www.ams.org/mathscinet-getitem?mr=#1}{#2}
} \providecommand{\href}[2]{#2}

\end{document}